%% file: cfrac.tex
\documentclass[onefignum,onetabnum]{siamart190516}
\usepackage{listings}

\usepackage{verbatimbox}

\usepackage{ifmtarg}
\usepackage{mathtools}
\usepackage{subcaption}
\usepackage{caption}
\usepackage{color}
\definecolor{dkgreen}{rgb}{0,0.6,0}
\definecolor{gray}{rgb}{0.5,0.5,0.5}
\definecolor{mauve}{rgb}{0.58,0,0.82}
\DeclarePairedDelimiter\ceil{\lceil}{\rceil}
\DeclarePairedDelimiter\floor{\lfloor}{\rfloor}

\DeclareMathOperator*{\argmin}{argmin}   
\makeatletter

    \newcommand\contFrac{\@ifstar{\@contFracStar}{\@contFracNoStar}}

    \def\singleContFrac#1#2{%
        \begin{array}{@{}c@{}}%
            \multicolumn{1}{c|}{#1}%
            \\%
            \hline%
            \multicolumn{1}{|c}{#2}%
        \end{array}%
    }

    \def\@contFracNoStar#1{%
        \mathchoice{
            \@contFracNoStarDisplay@#1//\@nil%
        }{
            \@contFracNoStarInline@#1//\@nil%
        }{
            \@contFracNoStarInline@#1//\@nil%
        }{
            \@contFracNoStarInline@#1//\@nil%
        }%
    }

    \def\@contFracNoStarDisplay@#1//#2\@nil{%
        \@ifmtarg{#2}{%
            #1%
        }{%
            #1+\cfrac{1}{\@contFracNoStarDisplay@#2\@nil}%
        }%
    }

        \def\@contFracNoStarInline@#1//#2\@nil{%
            \@ifmtarg{#2}{%
                #1%
            }{%
                #1 \@@contFracNoStarInline@@#2\@nil%
            }%
        }
        \def\@@contFracNoStarInline@@#1//#2\@nil{%
            \@ifmtarg{#2}{%
                + \singleContFrac{1}{#1}%
            }{%
                + \singleContFrac{1}{#1} \@@contFracNoStarInline@@#2\@nil%
            }%
        }

    \def\@contFracStar#1{%
        \mathchoice{
            \@contFracStarDisplay@#1////\@nil%
        }{
            \@contFracStarInline@#1//\@nil%
        }{
            \@contFracStarInline@#1//\@nil%
        }{
            \@contFracStarInline@#1//\@nil%
        }%
    }

    \def\@contFracStarDisplay@#1//#2//#3\@nil{%
        \@ifmtarg{#2}{%
            #1%
        }{%
            #1 + \cfrac{#2}{\@contFracStarDisplay@#3\@nil}%
        }%
    }

        \def\@contFracStarInline@#1//#2\@nil{%
            \@ifmtarg{#2}{%
                #1%
            }{%
                #1 \@@contFracStarInline@@#2\@nil%
            }%
        }
        \def\@@contFracStarInline@@#1//#2//#3\@nil{%
            \@ifmtarg{#3}{%
                + \singleContFrac{#1}{#2}%
            }{%
                + \singleContFrac{#1}{#2} \@@contFracStarInline@@#3\@nil%
            }%
        }
\makeatother

\input{ex_shared}

\ifpdf
\hypersetup{
  pdftitle={Numerical Continued Fraction Interpolation},
  pdfauthor={O. Salazar Celis}
}
\fi

\headers{Numerical Continued Fraction Interpolation}{O. Salazar Celis}

\title{Numerical Continued Fraction Interpolation
}

\author{Oliver Salazar Celis\thanks{Department of Mathematics and Computer Science. University of Antwerp, Antwerp, Belgium. 
  (\email{oliver.salazarcelis@uantwerpen.be}}), ING Belgium, Brussels, Belgium
  }

\begin{document}

\maketitle

{\centering \footnotesize  Dedicated to Silvia from her father.\par}

\begin{abstract}
We show that highly accurate approximations can often be obtained from constructing Thiele interpolating continued fractions by a Greedy selection of the interpolation points together with an early termination condition. The obtained results are comparable with the outcome from state-of-the-art rational interpolation techniques based on the barycentric form.
\end{abstract}


\section{Introduction}
An alternative title to this note could be \emph{The unreasonable effectiveness of Thiele continued fractions}. Indeed, it is known~\cite{CUYT1988,GM80} and further detailed in Section~\ref{num} 
that the construction of Thiele interpolating continued fractions can suffer from numerical instabilities. As argued in~\cite{GM80}, a careful selection of the ordering of the interpolation points is therefore needed.  
Here we propose a Greedy strategy, one which takes the next point where the error is maximal. It is shown in Theorem~\ref{thm:exist} that, at least from an existence point of view, this strategy is well motivated. We cannot prove that it is the overall \emph{best} strategy, but we rather give ample numerical evidence that the approach can lead to results that are competitive with other well-known approaches~\cite{Pachn2010AlgorithmsFP,Nakatsukasa2018TheAA,becker} for univariate rational interpolation and that the observed behavior is typical.

Once the existence problem is out of the way, we show how also best approximations can be obtained using the approach from~\cite{brasil} and how poles and zeros can be calculated directly from a generalised eigenvalue problem in Theorem~\eqref{thm:poles}. Combined, these cover most of the important tools required from a numerical rational interpolation method. 

\section{Thiele continued fraction rational interpolation}
\noindent Consider a sequence of distinct complex points $(x_i)_{i \in \mathbb{N}}$ together with function valuations $f(x_i)=f_i$  of the complex function $f(x)$.  
The continued fraction 
\begin{equation} \label{thielefrac}
\textstyle \contFrac*{   d_0//  x- x_{0} //  d_1 // x- x_{1} // d_2 // x- x_{2} // \dots }
\end{equation}
with $d_{i}= \varphi_{i}[x_0,\ldots,x_i] $ where
\begin{equation}
\left\{
\begin{aligned}
 \varphi_0[x_k] &= f_k \qquad &k \geq 0 \\
 \varphi_{i+1}[x_0,\ldots,x_i,x_k] &= \frac{x_k-x_i}{\varphi_{i}[x_0,\ldots,x_{i-1},x_k]-\varphi_{i}[x_0,\ldots,x_i]}\qquad &k > i 
\end{aligned}
\right. 
\end{equation}
generates rational interpolants when considering its successive convergents. The $n$th convergent 
\begin{equation} \label{eqn:thiele}
C_n(x) =  \textstyle \contFrac*{  \varphi_0[x_0]  //  x- x_{0} // \varphi_1[x_0,x_1] // x- x_{1} // \varphi_2[x_0,x_1,x_2] // x- x_{2} // \dots // x- x_{n-1}  // \varphi_n[x_0,\ldots,x_n]} 
\end{equation}
has the property that $C_n(x_i) = f_i$ for $i=0,\ldots,n$, provided that the inverse differences $\varphi_{i}[x_0,\ldots,x_i] \neq \infty$ exist and
none of the tails
\begin{equation}
\left\{
\begin{aligned}
 T_{i,n}(x)  &= \textstyle \contFrac*{  \varphi_i[x_0,\ldots,x_i]  //  x- x_{i} // \dots // x- x_{n-1} // \varphi_n[x_0,\ldots,x_n]} & \qquad 0 \leq i < n \\
  T_{n,n}(x) &= \varphi_n[x_0,\ldots,x_n] 
\end{aligned}
\right. 
\end{equation}
vanish at $x_{i-1}$ for $i > 0$. When $T_{i,n}(x_{i-1}) = 0$, then $x_{i-1}$ is \emph{unattainable}, meaning that it is a common zero of the numerator and denominator of $C_n(x)$ and in addition $\lim_{x \to x_{i-1}} C_{n}(x) \neq f_{i-1}$.

\section{Existence, sampling and numerical stability}
\subsection{Off-diagonal interpolants}
Recall~\cite[see~\S5.4, p.~110]{milne} that the numerator and denominator degrees of the convergents $C_n(x)$ are at most $\ceil*{n/2}$ and $\floor*{n/2}$ respectively.
 If desired, more off-diagonal interpolants can be obtained by considering $f(x) - p_\ell(x)$ instead of $f(x)$, with $p_\ell(x)$ an interpolating polynomial of degree at most $\ell$. The introduction of such a polynomial could be useful to \emph{dampen} functions $f(x)$ with a wide varying range. However, in what follows our interest is purely in Thiele continued fractions, hence we do not explore further potential benefits of non-diagonal interpolants.

\subsection{Unattainability}  While unattainability is undesirable, it is rarely encountered in practice except for pathological cases which are mostly related to symmetry of $f(x)$. Therefore we only mention that for $i = 1,\ldots,n$, it can be checked a posteriori if desired: either by checking for zeros of the tails $T_{i,n}(x)$ at $x_{i-1}$ or from direct evaluation of the convergents $C_n(x)$ at the interpolation points. The next theorem says that a Thiele continued fraction is irreducible, up to possible unattainable points.
\begin{theorem}
Provided that the convergents $C_n(x)$ exist, its only common factors are of the form $(x-x_i)$ with $x_i \in \{ x_0, \ldots,x_{n-1}\}$ an interpolation point. In such a case the $x_i$ are unattainable and are characterized by  vanishing of the tails $T_{i+1,n}(x_i)=0$.
\end{theorem}
\begin{proof}
Elaborate proofs can for instance be found in~\cite{werner,GM80}.
\end{proof}
\subsection{Existence} Non-existence of inverse differences on the other hand is more inconvenient  as it renders contributions from subsequent tails meaningless. The existence problem depends entirely on the ordering of the $(x_i)_{i \in \mathbb{N}}$. Specifically, it is required for two consecutive convergents $C_{i-1}(x)$ and $C_{i}(x)$ to be different ($i > 0 )$ in order for the inverse differences $\varphi_{i}[x_0,\ldots,x_i]$ to exist. 

\begin{theorem}\label{thm:exist}
If the points $(x_i)_{i \in \mathbb{N}}$ are ordered such that every two consecutive convergents of the continued fraction~\eqref{thielefrac} are different, then $\varphi_{i}[x_0,\ldots,x_i]  \neq \infty$.
\end{theorem}
\noindent The proof is detailed in Appendix~\ref{sec:proofs}, where it is shown that the inverse differences can essentially be interpreted as the ratio of two (linearized) residuals of successive convergents. 
This interpretation is used next to motivate a Greedy selection of the interpolation points. 

\subsection{Adaptive Greedy selection} \label{greedy} Given a finite sequence of points $(x_i)_{0\leq i \leq n}$, the condition of Theorem~\ref{thm:exist} where successive convergents should be different gives us a heuristic way to choose the next point $x_{i+1}$ in the construction of $C_{i+1}(x)$ with $0<i < n$. Given $C_{i}(x)$ and $(x_0,\ldots,x_i)$, reorder the remaining points $(x_{i+1}, \ldots,x_n)$ to determine $C_{i+1}(x)$ such that 
$|C_i(x_{i+1}) -f(x_{i+1})|$ is maximal. In this way, $C_n(x)$ is ultimate constructed in an adaptive Greedy way by choosing in each step the point where for $0<i<n$ the error between $C_{i}(x)$   and $f(x)$ is maximal. For the first point $x_0$ one can take a point where $|f(x_0)|$ is minimum. As such, at least one zero of $f(x)$ is accurately represented when present in the data. 

This Greedy selection ensures the existence of the inverse differences (see Appendix~\ref{sec:proofs}). A similar strategy is adopted by the AAA  approach~\cite{Nakatsukasa2018TheAA}, although there the choice is motivated by numerical considerations rather than existence. We cannot prove that the Greedy strategy is best for numerical stability in the context of Thiele continued fractions (discussed next), but the examples in this note suggest that the observed stable behaviour is typical.

\subsection{Numerical stability} \label{num}
It is known~\cite{CUYT1988} that the computation of the inverse differences and successive convergents in a Thiele interpolating continued fraction can suffer from numerical instabilities and worst case exponential loss of precision may occur in the calculation of the inverse differences~\cite{GM80}. Nevertheless, the backward evaluation of the continued fraction often leads to near machine precision (roughly \verb|2e-16|) magnitudes for the residuals $f(x_i) -C_n(x_i) $ in practice~\cite{CUYT1988}. 

We do not attempt an analysis like in~\cite{thron} because  any error analysis would necessarily depend on both the function $f(x)$ and the chosen interpolation points $x_i$.  For instance, a very strong stability result is recently obtained in~\cite{markovthiele} when restricting to Markov functions. Instead, we provide compelling numerical evidence in Sections~\ref{sec:newman} and~\ref{sec:other} that accurate approximations are often obtained, even on hard numerical problems, when the Greedy selection strategy of Section~\ref{greedy} is employed. 

\subsection{Early termination}  From the discussion in the previous section, it is important to obtain the Thiele interpolating continued fractions in as few steps as possible. Besides the greedy selection, we add a stopping criterion when constructing $C_n(x)$. If the maximum absolute error in the remaining points is below a prescribed tolerance, say \verb|tol=5e-15|, then we stop the recursion. The precise condition used is
$$
 \max |C_i(x_{k}) -f(x_{k})| < \text{tol} \times \max_{i<j\leq n}(|f(x_j)|) \qquad  \text{for } i < k \leq n.
$$It essentially tells us that within numerical tolerance, the function to be approximated is rational. This can also be understood from the deeper connection to the shape of the blocks in the Walsh table~\cite{Walsh1934}, a topic we do not address further.

\section{Example: interpolation of $|x|$ and $\sqrt{x}$}  \label{sec:newman}
 To illustrate that the Greedy selection strategy  works well in practice, we start with the example of interpolation of $f(x) = |x|$ in Newman~\cite{newman64} points for $x \in [-1,1]$.
This problem demonstrates the approximation power of rational functions as compared to polynomials which can only achieve $\mathcal{O}(n^{-1})$ accuracy at best. However, it has proven quite challenging numerically and has been used to assess the performance of several rational approximations schemes based on for instance the barycentric form~\cite{Pachn2010AlgorithmsFP,Nakatsukasa2018TheAA,Ion2020}. 

Newman approximations are rational interpolants in the $2n+1$ points
\begin{equation} \label{newpoints}
(-1, -\eta,\ldots,-\eta^{n-1}, 0, \eta^{n-1},\ldots, 1), \qquad \text{with } \eta = e^{-1/\sqrt{n}}.
\end{equation}
These points cluster around the origin, approaching it at an exponential rate as $n$ increases. 
It is shown in~\cite{zhou2004} that the asymptotic rate of convergence of rational interpolants to $|x|$ is root-exponential $\mathcal{O}( n^{-1/2} e^{-\sqrt{n}})$.

\begin{figure}[ht!]
     \centering
     \begin{subfigure}[b]{0.44\textwidth}
         \centering
         \includegraphics[width=\textwidth]{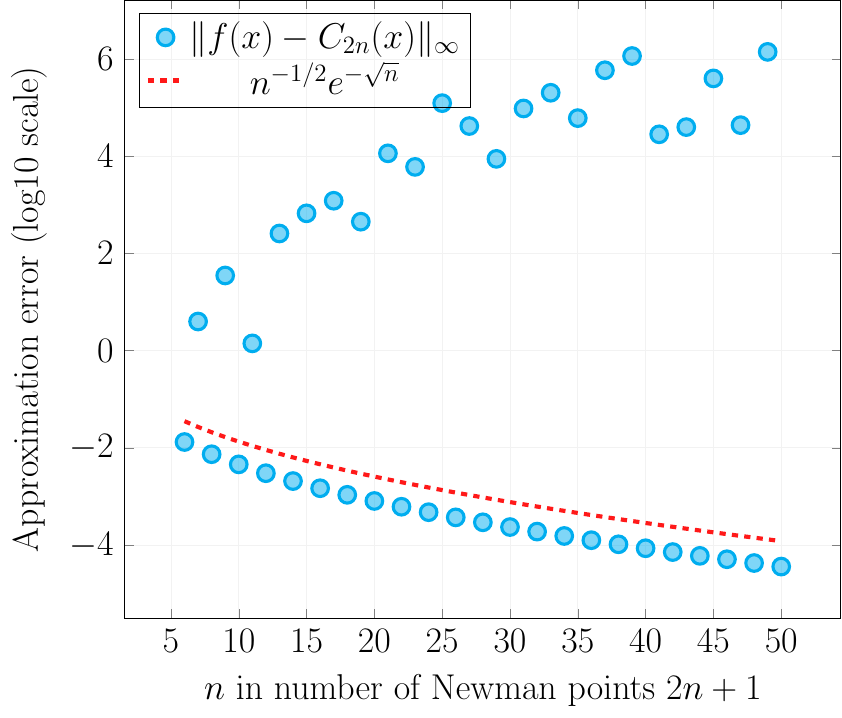}
         \caption{$\|f(x) - C_{2n}(x) \|_\infty$ on $x \in [-1,1] $}
         \label{fig:newmanapprox}
     \end{subfigure}
     \hfill
     \begin{subfigure}[b]{0.45\textwidth}
         \centering
         \includegraphics[width=\textwidth]{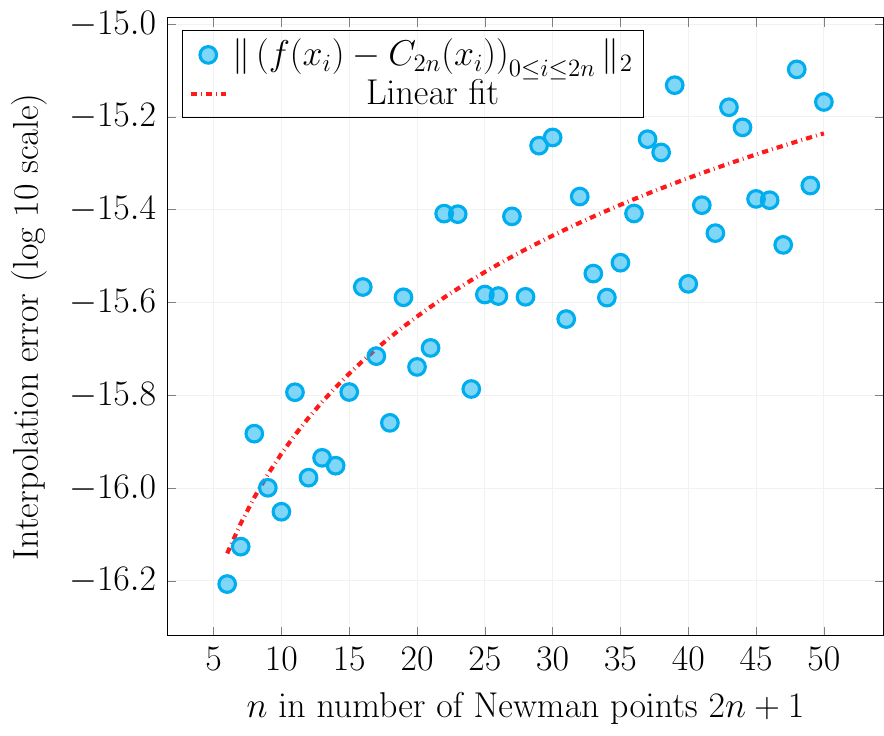}
         \caption{$\| \left( f(x_i)- C_{2n}(x_i) \right)_{0\leq i \leq 2n} \|_2$ }
         \label{fig:forwarderror}
     \end{subfigure}
        \caption{Adaptive Thiele interpolation of $|x|$ in Newman points for $n=6,\ldots,50$.  Left: maximum error on $[-1,1]$ in approximation of $|x|$ by Thiele interpolation. 
        Odd values of $n$ give approximations with poles in $[-1,1]$, while the even $n$ approximations are pole-free. The asymptotic bound~\cite{zhou2004} is superimposed (dotted line).
        Right: 2-norm of the interpolation error in the Newman points~\eqref{newpoints} from backward evaluation of the Thiele interpolant. The evaluations remain highly accurate around 15 decimal digits precision and errors grow linearly (dotted line) as $n$ increases. }
        \label{fig:absx}
\end{figure}

Figure~\ref{fig:absx} shows the results of Thiele interpolation in Newman points and $f(x) = |x|$ with $n$ up to $50$.  All obtained interpolants use all interpolation points in their construction, meaning that for $n=50$ we have $2n+1 = 101$ interpolation points and we construct $C_{2n}(x) = C_{100}(x)$. Mind that this is a challenging example because $C_{100}(x)$ is a rational function of degree $\ceil*{n} = 50$ and $\floor*{n}= 50$ in numerator and denominator respectively.
For instance the robust approach~\cite[see Fig.~5.6, p.~75]{Pachn2010AlgorithmsFP} on the same example already breaks down at $n \approx 11$ for reasons related to rounding errors. 
The Thiele interpolants on the other hand provide accurate results and follow the theoretical root-exponentional convergence. These results are in line with those observed in 
the more recent AAA approach~\cite[see Fig.~6.10, p.~1511]{Nakatsukasa2018TheAA}

\begin{figure}[h!]
     \centering
     \begin{subfigure}[b]{0.44\textwidth}
         \centering
         \includegraphics[width=\textwidth]{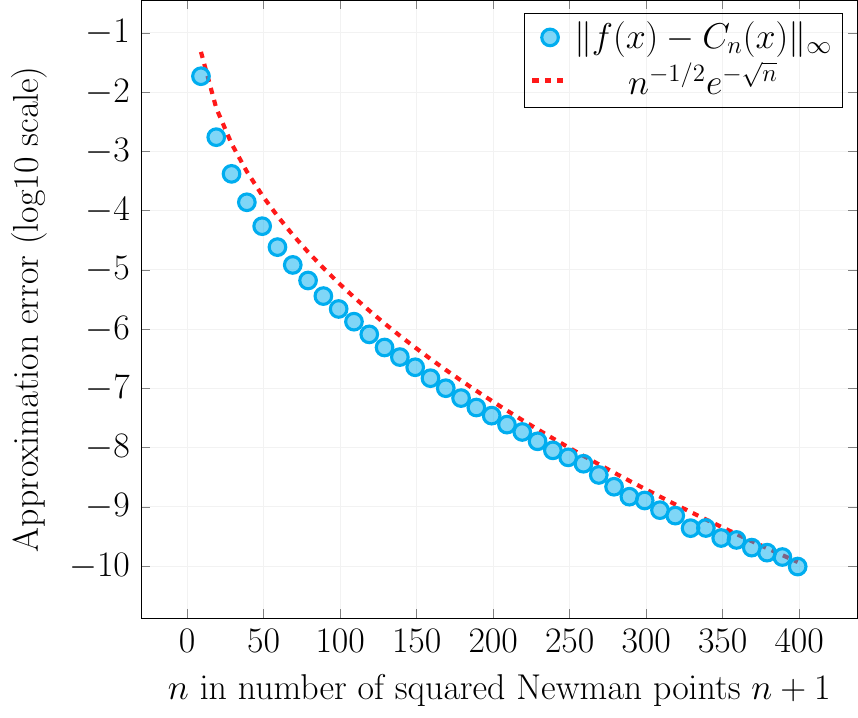}
         \caption{$\|f(x) - C_n(x) \|_\infty$ on $x \in [0,1] $}
         \label{fig:newmanapproxsqrt}
     \end{subfigure}
     \hfill
     \begin{subfigure}[b]{0.45\textwidth}
         \centering
         \includegraphics[width=\textwidth]{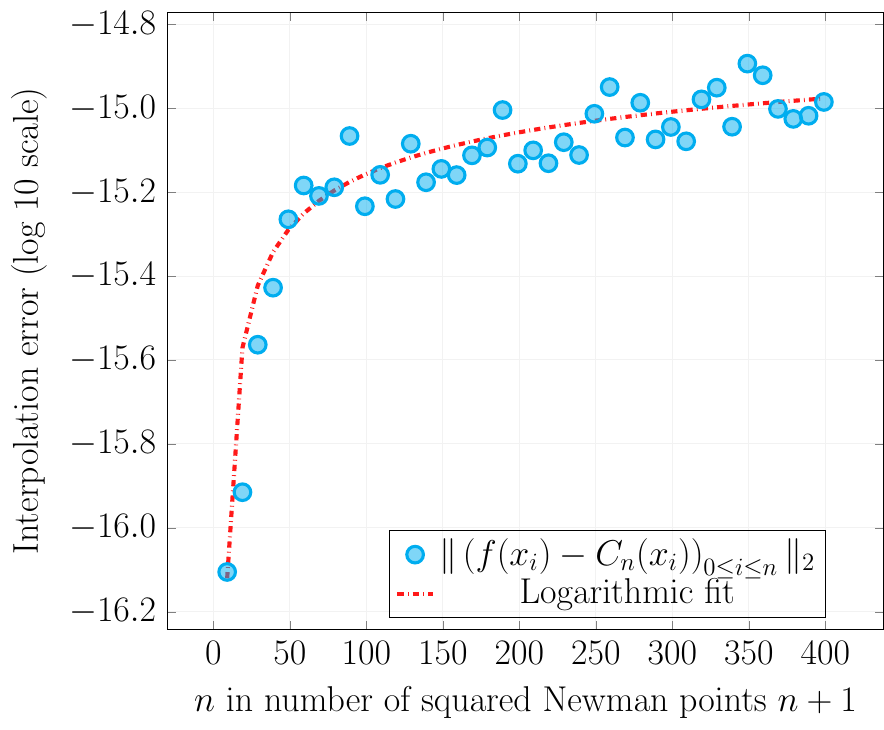}
         \caption{$\| \left( f(x_i)- C_{n}(x_i) \right)_{0\leq i \leq n} \|_2$ }
         \label{fig:forwarderror_sqrt}
     \end{subfigure}
        \caption{Adaptive Thiele interpolation of $\sqrt{x}$ in squared Newman points~\eqref{newpoints2} for $n=5,\ldots,400$.  Left: maximum error on $[0,1]$ in approximation of $\sqrt{x}$ by Thiele interpolation.  The asymptotic bound~\cite{zhou2004} is superimposed (dotted line).
        Right: 2-norm of the interpolation error in the chosen interpolation points from backward evaluation of the Thiele interpolant. The evaluations remain highly accurate around 15 decimal digits precision and errors grow logarithmically (dotted line) as $n$ increases.}
        \label{fig:sqrtx}
\end{figure}

An equivalent, but computationally easier, problem is to approximate $f(x) = \sqrt{x}$ on $x \in [0,1]$. For this interpolation problem we take the square of the Newmanpoints~\eqref{newpoints} leading to $n+1$ unique points in $[0,1]$:
\begin{equation} \label{newpoints2}
(0, \eta^{2(n-1)},\eta^{2(n-2)},\ldots, 1), \qquad \text{with } \eta = e^{-1/\sqrt{n}}.
\end{equation}
Note that these points cluster even closer near the origin than the original Newman points~\eqref{newpoints}.

Figure~\ref{fig:sqrtx} shows the results of Thiele interpolation in squared Newman points and $f(x) = \sqrt{x}$ with $n$ up to $400$.
The error in the interpolation points  seems to decay logarithmically with $n$. However, this is merely a consequence of the fact that for $n>60$ the interpolants are no longer of full degree. 
Beyond that point, not all given interpolation points are used in the continued fraction construction and their degree only rises slowly. In fact, for $n=400$ only $116$ points are used after which the early termination kicks in.
Nevertheless, the overall approximation quality improves in a root-exponential manner. We note that our implementation of the AAA method~\cite{Nakatsukasa2018TheAA} fails for this example when $n>50$ as it introduces real poles inside the interpolation interval.

\section{Best approximations}\label{sec:other}
While interpolation can deliver good approximations, interpolants are hardly ever \emph{best} approximations. For the remainder of this section, let $f(x) \in \mathcal{C}[a,b]$. If $C^*_n(x)$ is a rational function, then it is called a best approximation to $f(x)$ on $x \in [a,b]$ when
$$
 C^*_n(x) = \argmin_{C_n(x)} \max_{ x \in [a,b]} | f(x) - C_n(x) |. 
$$
Recent advances in best rational approximation have been made in~\cite{becker, NakatsukasaT20,brasil}.  We refer the interested reader to those references for further details on best approximations. 

Provided that $C^*_n(x)$ is of full degree, it is known that there exists a so-called \emph{alternant} set consisting of $n+2$ ordered nodes 
\begin{equation} \label{altern}
a \leq \tilde{x}_0 < \tilde{x}_1 <\cdots < \tilde{x}_{n+1} \leq b, 
\end{equation}
where $| f(x) - C^*_n(x) |$ attains its global extremum over all $x \in [a,b]$ with alternating signs:
$$
f(\tilde{x}_j) - C^*_n(\tilde{x}_j)  = (-1)^{\delta+j} \lambda, \qquad j=0,\ldots, n+1.
$$
Here $\delta \in  \{ 0,1\}$ and $\lambda = \max_{x \in [a,b]} | f(x) - C^*_n(x)|$. The alternating set is usually the starting point for the iterative Remez algorithm. From an interpolation point of view, the alternant set is not immediately useful because the construction of $C_n(x)$ requires only $n+1$ points and the value $\lambda$ is unknown. However, we can use an alternative approach as discussed next.

Due to continuity, the error $f(x) - C^*_n(x)$ must attain zero between each pair $(\tilde{x}_j,\tilde{x}_{j+1})$ of neighboring points of the alternant set. 
 This means that there must exist at least $n+1$ points  $x_0,\ldots,x_n$ 
 \begin{equation} \label{ref}
a \leq \tilde{x}_0 < x_0< \tilde{x}_1 < x_1 < \cdots < x_{n}< \tilde{x}_{n+1}  \leq b, 
\end{equation}
 such that $f(x_i) = C^*_n(x_i)$ for $i=0,\ldots,n$.
This idea is exploited in the BRASIL~\cite{brasil} algorithm which, given an initial guess of $n+1$ interpolation points, iteratively rescales the interval widths between successive interpolation points with the goal of equilibrating the local errors. In each step, a new set of interpolation points is determined. 

\begin{figure}[hbt!]
     \centering
     \begin{subfigure}[b]{0.46\textwidth}
         \centering
         \includegraphics[width=\textwidth]{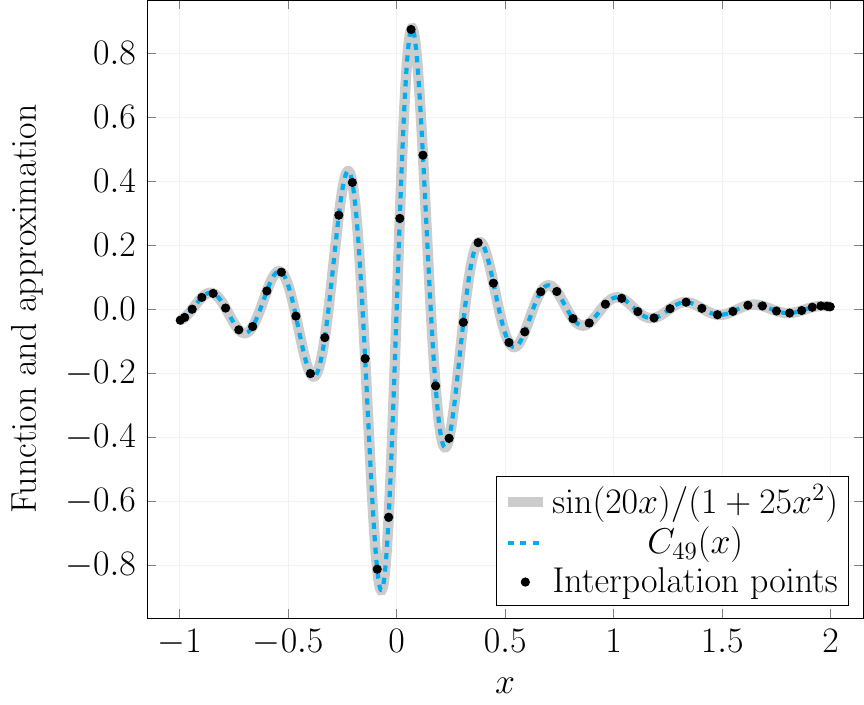}
         \caption{$f(x)$ and $C_{49}(x)$ on $x \in [-1,2] $}
         \label{fig:cosexp}
     \end{subfigure}
     \hfill
     \begin{subfigure}[b]{0.46\textwidth}
         \centering
         \includegraphics[width=\textwidth]{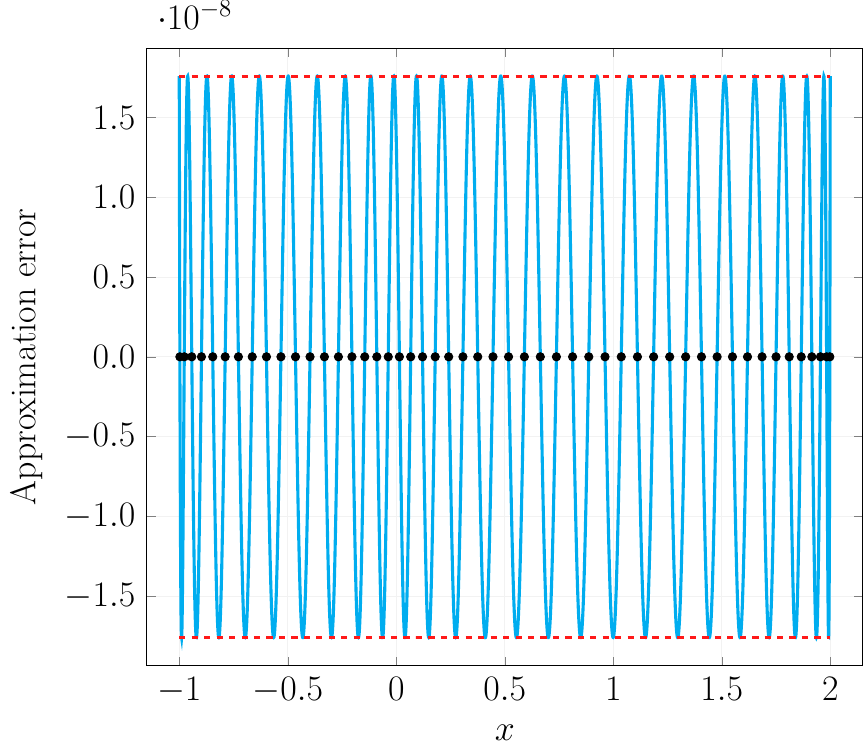}
         \caption{$f(x) - C_{49}(x) $ on $x \in [-1,2] $}
         \label{fig:cosexp_err}
     \end{subfigure}
     \begin{subfigure}[b]{0.46\textwidth}
         \centering
         \includegraphics[width=\textwidth]{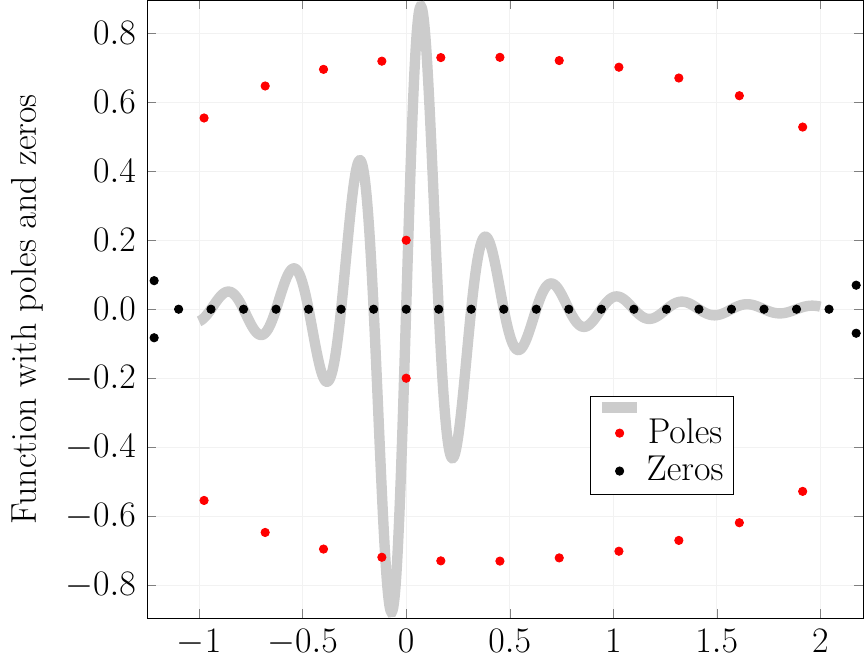}
         \caption{Poles and zeros of $C_{49}(x)$}
         \label{fig:poles_zeros}
     \end{subfigure}
        \caption{Best approximation of $f(x) =\sin(20x)/(1+25x^2)$ on $x \in [-1,2] $.  The interpolation points are shown with dots. Left: Best approximation $C_{49}(x)$ (dashed line), the function $f(x)$ is shown in gray. Right:  equioscillating error between $f(x)$ and $C_{49}(x)$. Bottom: Poles (red dots) and zeros (black dots) of $C_{49}(x)$. The function $f(x)$ is superimposed in grey. }
        \label{fig:sinratfit}
\end{figure}

Figure~\ref{fig:sinratfit} illustrates the best rational aproximation for $f(x) =\sin(20x)/(1+25x^2)$ on $x \in [-1,2] $. 
To obtain the initial interpolation points, a low accuracy continued fraction $C_{49}(x)$ is constructed from $100$ Chebychev points of the first kind between $-1$ and $2$. Then the BRASIL iteration is ran, which relocates the interpolation points in each step. The interpolation itself is done using the adaptive Thiele continued fractions approach rather than using the barycentric form as in~\cite{brasil}.  Other parameters, such as step size are put equal to $0.01$ and convergence acceleration was not implemented. The final error $f(x) - C_{49}(x)$ equioscillates $n+2 = 51$ times on  $x \in [-1,2] $ between the maximum error $\approx$ \verb|1.76e-08|.

Figure~\ref{fig:sinratfit} also shows the poles and zeros of the obtained continued fraction $C_{49}(x)$. 
One can clearly recognise that the two inner poles pick up the exact poles of $f(x)$ at  $\pm 0.2\mathrm{i}$. Also the zeros of $f(x)$ in the interval $[-1,2]$ are identified. The details of this calculation are given in Appendix~\ref{sec:poles}.

\begin{figure}[ht!]
     \centering
     \begin{subfigure}[b]{0.45\textwidth}
         \centering
         \includegraphics[width=\textwidth]{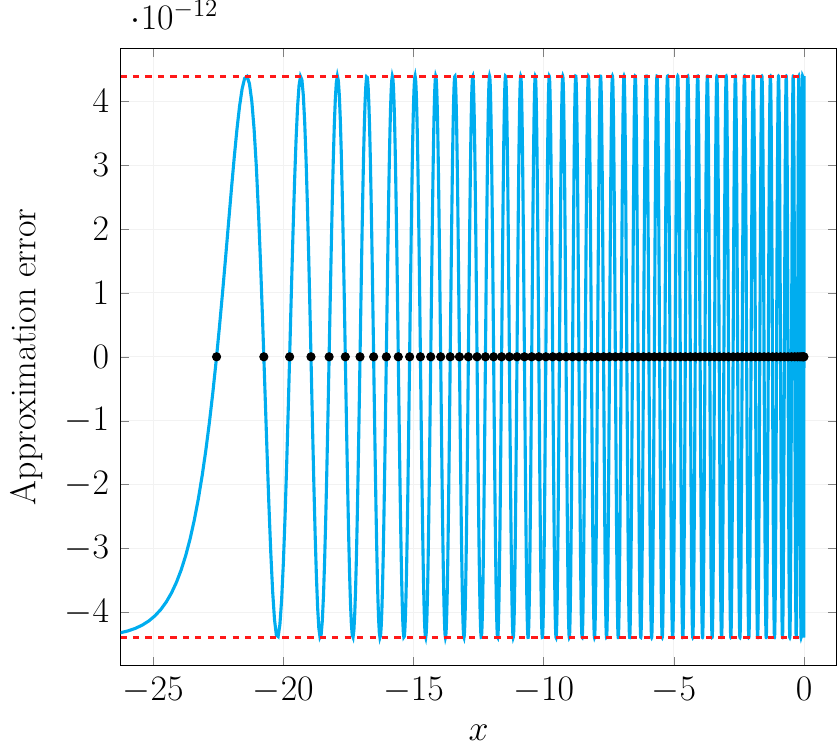}
         \caption{$\sqrt{x} - C^*_{80}(x) $ on $x \in [0,1] $}
         \label{fig:sqrt_log}
     \end{subfigure}
     \hfill
     \begin{subfigure}[b]{0.45\textwidth}
         \centering
         \includegraphics[width=\textwidth]{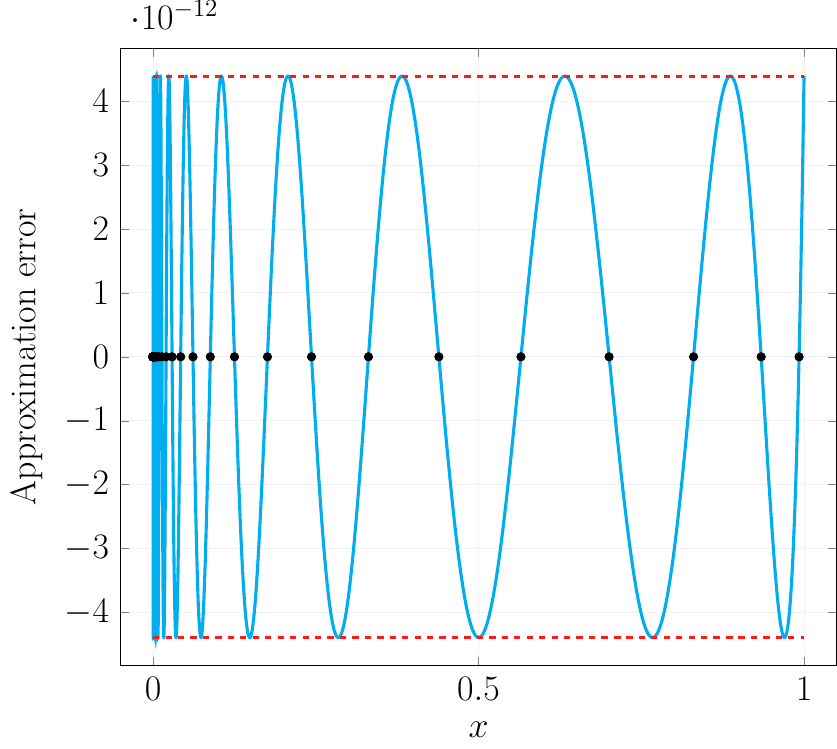}
         \caption{$\sqrt{x} - C^*_{80}(x)$ on $x \in [0,1] $}
         \label{fig:sqrt_err}
     \end{subfigure}
        \caption{Best approximation to $\sqrt{x}$ in $n+1=81$ points for $x \in [0,1]$.   The error from the final Thiele interpolant $C^*_{80}(x)$ equioscillates $n+2 = 82$ times between the maximum error (shown in dashed). The location of the interpolation points is shown with dots. Left: Log scale x-axis plot (base 10) of the equioscillating error to visualize the area close to the origin. Right: regular plot of the equioscillating error to visualize the area closer to $1$.} 
        \label{fig:sqrt}
\end{figure}

For the last illustration, we turn to the example of~\cite{varga}. In Figure~\ref{fig:sqrt}, the results for $f(x) = \sqrt{x}$ on $x \in [0,1]$ are shown. We found that a convenient starting point was obtained from a low(er) accuracy approximation based on $1000$ linearly spaced points between $0$ and $1$ raised to the power $6$ where the endpoints are removed. The iteration is started with $n+1=81$ points adaptively chosen from the aforementioned set  leading to the full degree interpolant  $C^*_{80}(x)$ with maximum leveled  error $\approx$ \verb|4.39e-12|. 
Compared to the interpolants of similar size (i.e.~constructed with the same number of interpolation points) in Figure~\ref{fig:newmanapproxsqrt}, the best approximation is roughly twice as accurate.

It is worth noting that one easily obtains an approximation of the same quality for $f(x)=|x|$ on $x \in [-1,1]$ by taking $C^*_{80}(x^2)$. Such an approximation would be equivalent to the one obtained in~\cite[Fig.~7.1]{becker} of degree $80$ in numerator and denominator. Mind that here the interpolation points range over more than 20 order of magnitudes and some of them are even below machine precision. A similar phenomenon is observed in~\cite[Fig.~3]{brasil}. To be fair, the direct best approximation of $f(x)=|x|$ on $x \in [-1,1]$ more often fails using the current approach, because the continued fraction representation struggles to maintain the exact symmetry.

\newpage
\section{Conclusions}
We have shown that highly accurate approximations can often be obtained with Thiele continued fractions when incorporating a Greedy selection of the interpolation points together with an early termination condition. The obtained results are comparable with state-of-the art rational interpolation techniques based on the barycentric form such as AAA~\cite{Nakatsukasa2018TheAA} and can be used for rational minimax approximation~\cite{brasil}. Also the poles and zeros are relatively easy to obtain, making the approach attractive for practical purposes.

The main advantage of the Thiele continued fraction is the simplicity of its construction. It does not require specialized linear algebra  implementations such as SVD, but rather relies on an elementary recursion (see Appendix~\ref{sec:rcode} for an implementaion). Of course accumulation of rounding errors does occur, hence inevitably one can find examples where the approach will break down. The possibility of  numerical breakdown is applicable for all rational interpolation approaches, including AAA~\cite[see \S6 p.3172--3173]{NakatsukasaT20}. 


\appendix
\section{Proof of Theorem~\ref{thm:exist}} \label{sec:proofs}

For the proof of Theorem~\ref{thm:exist} we first relate the construction of the inverse difference to the (linearized) residuals of the convergents $C_n(x)$ of the continued fraction~\eqref{thielefrac}.
To that end, let $C_n(x) = A_n(x)/B_n(x)$ where the $n$th numerator $A_n(x)$ and $n$th denominator $B_n(x)$ satisfy the recurrence relation~\cite{wallis}
\begin{equation} \label{eqn:3term}
\left(
\begin{aligned}
A_n(x) \\
B_n(x)
\end{aligned}
\right)
=\left(
\begin{aligned}
\varphi_{n}[x_0,\ldots,x_n]  A_{n-1}(x) +(x-x_{n-1})A_{n-2}(x)\\
\varphi_{n}[x_0,\ldots,x_n]  B_{n-1}(x) +(x-x_{n-1})B_{n-2}(x)
\end{aligned}
\right),
\end{equation}
with
$$
\left\{
\begin{aligned}
A_{-2}(x) = 0, & \qquad B_{-2}(x) = 1 \\
A_{-1}(x) = 1, & \qquad B_{-1}(x) = 0 \\
A_{0}(x) = \varphi_{0}[x_0] = f(x_0), &\qquad  B_{0}(x) = 1
\end{aligned}
\right. .
$$
\begin{theorem} \label{thm:res} If the points $(x_i)_{i \in \mathbb{N}}$ are ordered such that every two consecutive convergents of the continued fraction~\eqref{thielefrac} are different then
for $i = -1,0,1,\ldots $
$$
\left\{
\begin{aligned}
 \varphi_{0}[x_k] = - \frac{ R_{-2}(x_k)}{R_{-1}(x_k) }, \qquad k > i =-1 \\
 \varphi_{i+1}[x_0,\ldots,x_i,x_k] = - (x_k-x_i)\frac{ R_{i-1}(x_k)}{R_i(x_k) }, \qquad k > i \geq 0
\end{aligned}
\right. ,
$$
where the linearized residuals $R_i(x)$ are defined as
$$R_i(x) = f(x)B_i(x) - A_i(x).$$
\end{theorem}
\begin{proof}
The proof is by induction. For $i=-1$
$$
 - \frac{ R_{-2}(x_k)}{R_{-1}(x_k) }= -\frac{f(x_k)B_{-2}(x_k) -A_{-2}(x_k)}{f(x_k)B_{-1}(x_k)-A_{-1}(x_k)} = -\frac{f(x_k)}{-1} = \varphi_{0}[x_k] .
$$
And for $i=0$
$$
\begin{aligned}
 - (x_k-x_0)\frac{ R_{-1}(x_k)}{R_{0}(x_k) } &= -(x_k-x_0)\frac{f(x_k)B_{-1}(x_k) -A_{-1}(x_k)}{f(x_k)B_{0}(x_k)-A_{0}(x_k)}  \\
 & = -(x_k-x_0)\frac{-1}{f(x_k)-f(x_0)} =  \varphi_{1}[x_0,x_k] .
\end{aligned}
$$
Assume that the hypothesis  $\varphi_{i}[x_0,\ldots,x_{i-1},x_k] = - (x_k-x_{i-1}){ R_{i-2}(x_k)}/{R_{i-1}(x_k) }$ holds.  First note that from application of~\eqref{eqn:3term} we have
$$
R_i(x_k) = (x_k-x_{i-1}) R_{i-2}(x_k) + \varphi_{i}[x_0,\ldots,x_i]R_{i-1}(x_k).
$$
Hence
$$
\begin{aligned}
 - (x_k-x_i)\frac{ R_{i-1}(x_k)}{R_{i}(x_k) } = \frac{- (x_k-x_i) R_{i-1}(x_k)}{(x_k-x_{i-1}) R_{i-2}(x_k) + \varphi_{i}[x_0,\ldots,x_i]R_{i-1}(x_k)}
\end{aligned}
$$
Provided that $R_{i-1}(x_k) \neq 0$ and application of the induction hypothesis gives
$$
 - (x_k-x_i)\frac{ R_{i-1}(x_k)}{R_{i}(x_k) }  = \frac{-(x_k-x_i)}{-\varphi_{i}[x_0,\ldots,x_{i-1},x_k]+\varphi_{i}[x_0,\ldots,x_i]} =  \varphi_{i+1}[x_0,\ldots,x_{i},x_k].
$$
For the last step, the presumed ordering of the points is important. If two consecutive convergents $C_{i-1}(x)$ and $C_{i}(x)$ are different, then the polynomial 
\begin{equation} \label{eqn:nonzero}
A_{i-1}(x)B_{i}(x) -A_{i}(x)B_{i-1}(x) \not\equiv 0 
\end{equation}
is of degree at most $\max\{ \ceil{(i-1)/2}+\floor{i/2}, \ceil{i/2}+\floor{(i-1)/2} \} \leq i$ which already vanishes at $i$ points $(x_j)_{0 \leq j \leq i-1}$
$$
 R_{i-1}(x_j)B_{i}(x_j) - R_{i}(x_j)B_{i-1}(x_j) = 0, \qquad j=0,\ldots,i-1.
$$
Due to the condition~\eqref{eqn:nonzero} it cannot be that both $R_{i-1}(x)$ and $R_i(x)$ simultaneously vanish at other points. 
\end{proof}

\begin{proof}[Proof of Theorem~\ref{thm:exist}]
From Theorem~\ref{thm:res}, we know that the inverse differences can be expressed as
$$
 \varphi_{i+1}[x_0,\ldots,x_i,x_{i+1}] = - (x_{i+1}-x_i) \frac{ R_{i-1}(x_{i+1})}{R_i(x_{i+1}) }, \qquad i\geq0
 $$
 and the presumed ordering of the $(x_i)_{i \in \mathbb{N}}$ ensures that the ratio ${ R_{i-1}(x_{i+1})}/{R_i(x_{i+1})}$ is well-defined because its numerator and denominator cannot vanish simultaneously.
 In addition, $\varphi_{i+1}[x_0,\ldots,x_i,x_{i+1}] =\infty$ would imply  $R_i(x_{i+1})=0$ which contradicts that consecutive convergents $C_{i+1}(x)$ and $C_{i}(x)$ are different.
\end{proof}

\section{Poles, zeros and residues} \label{sec:poles}
The poles and zeros of a Thiele continued fraction can be extracted directly from application of a (generalised) eigenvalue problem. 
The key ingredient to that end is the \emph{continuant}~\cite{determinants} representation of the partial numerator and denominator.  
Even though continuants are well-known in the continued fraction literature, to the best of our knowledge their use to extract poles and zeros has not been exploited before.
\subsection{Poles}
\begin{theorem} \label{thm:poles}
Given distinct points $x_i$, finite $\varphi_i[x_0,\ldots,x_i]\neq \infty$ $(i=0,\ldots,n)$ 
and its associated Thiele continued fraction  
$$
C_n(x) =  \textstyle \contFrac*{  \varphi_0[x_0]  //  x- x_{0} // \varphi_1[x_0,x_1] // x- x_{1} // \varphi_2[x_0,x_1,x_2] // x- x_{2} // \dots // x- x_{n-1}  // \varphi_n[x_0,\ldots,x_n]}. 
$$
\begin{equation} \label{highestdegree}
\text{In case that } n \text{ is odd, assume:} \qquad \sum_{i=0, i \text{odd}}^{n} \varphi_i[x_0,\ldots,x_i] \neq 0.
\end{equation}
Let the $n\times n$ matrices ${\bf{D}}$ and $\bf{C}$ respectively be defined as
$$
{\bf{D}} = 
\begin{pmatrix}
\varphi_1[x_0,x_1] & -x_1 & 0 & \cdots & 0 \\
-1 & \varphi_2[x_0,x_1,x_2] &  -x_2  &  \ddots& \vdots  \\
0  & -1 & \varphi_3[x_0,x_1,x_2,x_3] & \ddots & 0 \\
\vdots & \ddots  & \ddots & \ddots & -x_{n-1}\\
0 & \cdots & 0 & -1 & \varphi_n[x_0,\ldots,x_n] 
\end{pmatrix},
$$
$$
{\bf{C}}=
\begin{pmatrix}
0 & -1 & 0 & \cdots & 0 \\
0 &  0&  -1  &  \ddots& \vdots  \\
0  & 0 & 0 & \ddots & 0 \\
\vdots & \ddots  & \ddots & \ddots & -1\\
0 & \cdots & 0 & 0& 0
\end{pmatrix}
$$
then the matrix pencil ${\bf{D}}-\lambda {\bf{C}}$ is regular and its finite eigenvalues coincide with the $\floor*{n/2}$ poles of $C_n(x)$ counting multiplicities.
\end{theorem}
\begin{proof}
Using the three term recurrence~\eqref{eqn:3term}, it is well-known~\cite{determinants} that one can write the denominator $B_n(x)$ of $C_n(x)$ as the determinant of an $n \times n$ tridiagonal matrix, which is also called a \emph{continuant}
\begin{equation} \label{eqn:continuant}
B_n(x) = 
\begin{vmatrix}
\varphi_1[x_0,x_1] & x-x_1 & 0 & \cdots & 0 \\
-1 & \varphi_2[x_0,x_1,x_2] &  x-x_2  &  \ddots& \vdots  \\
0  & -1  & \varphi_3[x_0,x_1,x_2,x_3]  & \ddots & 0 \\
\vdots & \ddots  & \ddots & \ddots & x-x_{n-1}\\
0 & \cdots & 0 & -1 & \varphi_n[x_0,\ldots,x_n] 
\end{vmatrix}.
\end{equation}
Clearly~\eqref{eqn:continuant} can be written in the form  $\text{det}({\bf{D}}-x{\bf{C}})$ where the $n\times n$ matrices $\bf{D}$ and $\bf{C}$ are as defined above.
Under the given assumptions, $B_n(x)$ is of exact degree $\floor*{n/2}$. Hence, $B_n(x) = \text{det}({\bf{D}}-x{\bf{C}}) \not\equiv 0$ so that the matrix pencil ${\bf{D}}-\lambda {\bf{C}}$ is regular.
 After all, in case that $n$ is even, then $B_n(x)$ is monic  of exact degree $n/2$; regardless of the  values $\varphi_i[x_0,\ldots,x_i]\neq \infty$ $(i=0,\ldots,n)$. 
 In case that $n$ is odd, then $B_n(x)$ is of exact degree $(n-1)/2=\floor*{n/2}$ with highest degree coefficient equal to the sum in~\eqref{highestdegree}.
 
 Since the characteristic polynomial $\text{det}({\bf{D}}-\lambda{\bf{C}})= B_n(\lambda)$, the poles of $C_n(x)$ coincide with the non-zero eigenvalues of ${\bf{D}}-\lambda {\bf{C}}$. 
 There are exactly $\floor*{n/2}$ such eigenvalues, because $B_n(\lambda)$ is of exact degree $\floor*{n/2}$.
\end{proof}

{\bf{Remark.1}} The condition~\eqref{highestdegree} is seldom of practical importance. It is more of theoretical interest to prevent the trivial case $B_n(x) \equiv 0$ and the occurrence of additional eigenvalues at infinity when $B_n(x)$ is not of exact degree $\floor*{n/2}$. 

{\bf{Remark 2.}} In practice, the poles of $C_n(x)$ can thus be found by solving the generalised eigenvalue problem
$$
{\bf{D}v}= \lambda {\bf{C}v}.
$$
The eigenvalues at infinity can be discarded. There are at least $\ceil*{n/2}$ infinite eigenvalues (also see Remark 1).
\subsection{Zeros}
In an entirely similar fashion, also the roots of $C_n(x)$ can be determined from a generalised eigenvalue problem. Based on the continuant representation of the numerator of $C_n(x)$, one simply replaces the $n\times n$ matrices ${\bf{D}}$ and $\bf{C}$ by the $(n+1) \times (n+1)$ matrices, where ${\bf{D}}$ is now
$$
{\bf{D}} = 
\begin{pmatrix}
\varphi_0[x_0] & -x_0 & 0 & \cdots & 0 \\
-1 & \varphi_1[x_0,x_1] &  -x_1  &  \ddots& \vdots  \\
0  & -1 & \varphi_2[x_0,x_1,x_2] & \ddots & 0 \\
\vdots & \ddots  & \ddots & \ddots & -x_{n-1}\\
0 & \cdots & 0 & -1 & \varphi_{n}[x_0,\ldots,x_{n}] 
\end{pmatrix}.
$$
and swaps every occurrence of odd by even in condition~\eqref{highestdegree} of Theorem~\eqref{thm:poles}. Recall that $C_n(x)$ has at most $\ceil*{n/2}$ zeros counting multiplicities.
\subsection{Residues}
Once the poles are determined, the residues can be found in a similar fashion as in~\cite{Nakatsukasa2018TheAA}. Essentially, one applies a trapezoidal rule approximation to
the contour integral of $C_n(x)$ over a small circle around each pole in the complex plane.

\section{R Code} \label{sec:rcode}
All codes have been developed in the open source R language~\cite{rmanual}. The basic code for performing Thiele interpolation together with an example is given below. 
\begin{verbnobox}[\small]
evalcfrac <- function(aa,zz,x){
  #backward evaluation of continued fraction aa1 + (x-zz1)/aa2+...
  j = max(which(is.finite(aa))); res=rep(0,length(x))
  if(j>1){for(i in j:2){res=(x-zz[i-1])/(aa[i]+res)}}
  return(aa[1]+res)}

cfrac_interpolate <-function(xx,ff,tol=5e-15,NN=0){
  if(NN==0){NN=length(xx);} #use tol criteria to stop
  aa=rep(NA,NN); zz=rep(NA,NN); rr=rep(NA,NN)
  for(k in 1:NN){  #main loop
    if(k==1){ #init
      rr=ff #inverse differences
      i=which.min(abs(ff)) #smallest value
    }else{ i=which.max(abs(evalcfrac(aa,zz,xx)-ff)) #adaptive choice
    rr=(xx-zz[k-1])/(rr-aa[k-1])} #inverse differences
    aa[k]=rr[i];zz[k]=xx[i] #store cfrac coef
    ff=ff[xx!=xx[i]];rr=rr[xx!=xx[i]];xx=xx[xx!=xx[i]]; #reduce data
    if(k<NN){
      if(max(abs(evalcfrac(aa,zz,xx)-ff))<tol*max(abs(ff))){
        print(paste("target precision reached earlier at n=",k))
        break();}}
  };return(list( "a"=aa[is.finite(aa)],"z"=zz[is.finite(zz)]))}

#example
f <- function(x) cos(exp(x))
x=seq(-1,1,length.out=100); y=f(x)
cfr=cfrac_interpolate(x,y)
plot(x,abs(f(x)-evalcfrac(cfr$a,cfr$z,x)))
\end{verbnobox}

The basic implementation of the BRASIL~\cite{brasil} iteration is given below. 

\begin{verbnobox}[\small]
findExtrema <- function(xk,a,b,cfr,f){
  xk = sort(unique(c(a,xk,b))) #add endpoints
  sk = numeric(0);resids = numeric(0);residabs = numeric(0)
  for(i in 2:length(xk)){ ## get extrema
    obj_func <- function(x) (f(x)-evalcfrac(cfr$a,cfr$z,x))
    opt = optimize(f = function(x) abs(obj_func(x)), 
                lower = xk[i-1], upper = xk[i],maximum = TRUE, tol = 1e-30)
    resids = c(resids, f(opt$maximum)-evalcfrac(cfr$a,cfr$z,opt$maximum))
    sk = c(sk, opt$maximum)
  };res <- list("sk"= sk,"resids"= resids,"residabs" = abs(resids)) }
  
brasil <- function(xk,a,b,cfr,f,smax=1e-1,tt=1e-1, tol=1e-3){
  stop = FALSE
  while(!stop){
    oldx = cfr$z[is.finite(cfr$z)]
    res = findExtrema(oldx,a=a,b=b,cfr,f)
    eps = max(res$residabs)/min(res$residabs)-1 
    if(eps<tol){ stop =TRUE
    }else{
      h = mean(res$residabs) #mean err
      g = max(abs(res$residabs-h)) #max deviation from mean err
      gk = (res$residabs-h)/g
      s  = min(c(smax, tt*g/h)) #stepzise
      ck = (1-s)^gk
      xkk = unique(sort(c(a,oldx,b)))
      lk = ck*diff(xkk) #rescaled length
      w = sum(sort(lk)) #normalisation
      newx = rep(NA, length(oldx))
      for(i in 1:length(newx)){
        newx[i] = (a*w + (b-a)*sum(lk[1:i]))/w
      };cfr = cfrac_interpolate(newx,f(newx),tol=0)
    }
  };return(cfr)}  
#example sqrt
f <- function(x) sqrt((x))
x=unique(seq(0,1,length.out=1000)^6); x = sort(x)[2:(length(x)-1)]
cfr = cfrac_interpolate(x,f(x),tol=5e-15, NN=81) #init
cfrstar = brasil(cfr$z,0,1,cfr,f,smax=1e-1,tt=1e-1, tol=2e-4)
plot(x,abs(f(x)- evalcfrac(cfrstar$a, cfrstar$z,x)),type="l")
\end{verbnobox}

Poles, zeros and residues are calculated as follows.

\begin{verbnobox}[\small]
prz_cfrac <-function(aa,zz){ # compute poles, residues, zeros
  j = max(which(is.finite(aa)));
  C = diag(0,j-1);diag(C[,-1]) = -1 ;
  Dp =diag(aa[2:j]);diag(Dp[,-1])=-zz[2:(j-1)];diag(Dp[-1,])=-1
  poleig = geigen(Dp,C,only.values=TRUE) 
  pol = poleig$values[is.finite(poleig$values)] #poles
  dz = 1e-5*exp(2i*pi*(1:4)/4)
  pz = outer(pol,dz,'+'); m = length(pol)
  res = matrix(data=NA, nrow = m, ncol=1) #residues
  for (k in 1:m){res[k]=evalcfrac(aa,zz, pz[k,])
  C = diag(0,j); diag(C[,-1]) = -1 ;
  Dz =diag(aa[1:j]); diag(Dz[,-1])=-zz[1:(j-1)];diag(Dz[-1,])=-1
  zereig = geigen(Dz,C, only.values=TRUE)
  zer = zereig$values[is.finite(zereig$values)] #zeros
  res <- list("pol"= pol, "zer" = zer, "res" = res)
}
\end{verbnobox}

\bibliographystyle{siamplain}
\bibliography{cfrac_bib2}
\end{document}

%% file: ex_shared.tex
\usepackage{lipsum}
\usepackage{amsfonts}
\usepackage{graphicx}
\usepackage{epstopdf}
\usepackage{algorithmic}
\ifpdf
  \DeclareGraphicsExtensions{.eps,.pdf,.png,.jpg}
\else
  \DeclareGraphicsExtensions{.eps}
\fi

\newsiamremark{remark}{Remark}
\newsiamremark{hypothesis}{Hypothesis}
\crefname{hypothesis}{Hypothesis}{Hypotheses}
\newsiamthm{claim}{Claim}

\usepackage{amsopn}
